\documentstyle[12pt]{article}
\title{Duality Theorems for Infinite Braided Hopf Algebras }
\author{
Shouchuan Zhang $^{a,~b}$, \ \  Yanying Han $^a$ \ \  \ \\$a$. Department  of Mathematics,
Hunan University\\ Changsha  410082, \
 P.R. China. E-mail:Zhangsc@hunu.edu.cn \\
$b$. Department of Mathematics, University of Queensland\\
Brisbane 4072, Australia. E-mail:yzz@maths.uq.edu.au  \\  }
\date{}
\begin{document}
\newtheorem{Theorem}{\quad Theorem}[section]
\newtheorem{Proposition}[Theorem]{\quad Proposition}
\newtheorem{Definition}[Theorem]{\quad Definition}
\newtheorem{Corollary}[Theorem]{\quad Corollary}
\newtheorem{Lemma}[Theorem]{\quad Lemma}
\newtheorem{Example}[Theorem]{\quad Example}
\maketitle \addtocounter{section}{-1}

 \begin {abstract} Let $H$ be an infinite-dimensional braided  Hopf algebra
 and assume that the braiding is symmetric on $H$ and its quasi-dual $H^d$.
We prove the Blattner-Montgomery  duality theorem, namely we prove
$$   (R \# H)\# H^{d}   \cong R \otimes (H \#  H^{d} ) \hbox {\ \ \ as algebras in braided tensor category  } {\cal C}.$$ In particular,
we present two duality theorems for infinite braided Hopf algebras in the Yetter-Drinfeld module category.

 \end {abstract}

\section {Introduction} The duality theorems play an important role in actions of Hopf algebras (see \cite {Mo93}).
In \cite {BM85} and \cite {Mo93},  Blattner and Montgomery
proved the following  duality theorem for an ordinary Hopf algebra $H$ and some Hopf subalgebra
$U$ of $H^\circ$ :
\begin {eqnarray*} \label {e1}
(R \# H)\# U   \cong R \otimes (H \# U)    \hbox { \ \ \ as algebras,
\  }
 \end {eqnarray*}
 where $R$ is a $U$-comodule algebra. The dual theorems for co-Frobenius Hopf algebra $H$,
$$   (R \# H)\# H^{* rat}   \cong M_H^f (R) \ \hbox { and }
(R \# H^{* rat })\# H   \cong M_H^f (R)\hbox { \ \ \  as  } k
\hbox {-algebras  } $$
were proved in \cite {DZ99} (see \cite [Corollary 6.5.6 and Theorem 6.5.11  ]
{DNR01}). On the other hand, braided
Hopf algebras have attracted much attention in both mathematics and mathematical physics
(see  \cite {AS02}\cite {Ke99} \cite {Ma95b}).
 One of the authors in \cite {Zh03} generalized the duality theorem
to the braided case, i.e., for a finite  Hopf algebra $H$
with  $C_{H,H} = C_{H,H}^{-1},$
  \begin {eqnarray*}
(R \# H)\# H^{\hat *}   \cong R \otimes (H \bar
  \otimes H^{\hat *} )  \hbox { \ \ \ as
algebras in } {\cal C}.
 \end {eqnarray*} The Blattner-Montgomery  duality theorem was also generalized to Hopf algebras over commutative rings
\cite {ATL01}.

In this paper we generalize the above results   to infinite braided Hopf
algebras. A braided Hopf algebra is called an infinite braided Hopf algebra if it has no
left duals (See \cite {Ta99}). An important example of
infinite braided Hopf algebras is the universal  enveloping algebra of a Lie
superalgebra. So  the dality theorems  for infinite braided Hopf
algebras should
have important applications in both mathematics and mathematical physics.

This paper is organized as follows.
In section 1, since it is possible that $Hom (H, I)$ is not an object in ${\cal C}$ for braided Hopf algebra $H$,
 we introduce  quasi-dual $H^d$ of  $H$,  and
prove the duality theorem in a braided tensor category  ${\cal C}$, i.e.
$(R \# H)\# H^{d}   \cong R \otimes (H \#  H^{d})$  as algebras  in  ${\cal C}$.
In section 2, we concentrate on the Yetter-Drinfeld module category $^B_B {\cal YD}$, and show
$$   (R \# H)\# U   \cong R \otimes (H \#  U )  \hbox {\ \ \  and  \ \ \ }
(R \# U)\# V   \cong R \otimes (U \#   V ) $$
 as algebras in $^B_B {\cal YD}.$ Here $U$ and $V$ are certain  braided Hopf subalgebras of $H^\circ$ and $H$, respectively.

\section {Duality theorem for braided Hopf algebras}

In this section, we obtain the duality theorem
 for braided  Hopf algebras living in the braided tensor category ${\cal C}.$

Let  $({\cal C}, \otimes, I,  C )$ be   a braided tensor category,
where $I$ is the identity object and $C$ is the braiding. we  write $W \otimes f $ for
$id _W \otimes f$ and $f \otimes W$ for  $f \otimes id _W$. The proofs in this section are
 very similar to those  for  the corresponding results in \cite [Chapter 9] {Mo93}
  and  \cite [Chapter 7] {Zh99}, so we only give the sketch to the proofs. In particular,
 there are the proofs in \cite {Zh99} by using braiding
 diagrams.

\begin {Definition} \label {1.1}
Let $(H, m, \eta , \Delta , \epsilon )$ in ${\cal C}$ be a braided Hopf algebra. If there is a braided Hopf algebra  $H^d$
and a morphism $<,>$  from $H^d\otimes H$ to $I$ such that
  \begin {eqnarray*}
&{\ }&( <, > \otimes <, >)(H^d \otimes C \otimes H)(\Delta _{H^d}\otimes H \otimes H) = <, >( H ^d\otimes m _H),\\
&{\ }&\epsilon _{H^d} = <, > (H^d \otimes \eta  _H), \\
&{\ }&<, >(m_{H^d} \otimes H)=( <, > \otimes <, >)(H^d \otimes C \otimes H)( H \otimes H \otimes \Delta ), \\
&{\ }&<, >(\eta _{H^d}\otimes H) = \epsilon _H , \\
 &{\ }&<, > (S_{H^d} \otimes H) = <,> (H^d\otimes S_H),\\
\end {eqnarray*}
 then $H^d$ is called a quasi-dual of $H$ under $<,>$.
\end {Definition}

\begin {Lemma} \label {1.2} Let $H^d$ be a quasi-dual of $H$ under $<,>$ and
$C_{H,H} = C_{H,H}^{-1}$.
Assume that $(H \otimes <,>)(C_{H^d, H} \otimes H)=
(H \otimes <,>)(C_{H, H^d} ^{-1} \otimes H)$ implies $C_{H^d, H}=
C_{H, H^d} ^{-1} $,  and $(H \otimes <,>)(C_{H^d, H ^d} \otimes H)=
(H^d \otimes <,>)(C_{H^d, H^d} ^{-1} \otimes H)$ implies $C_{H^d, H^d}=
C_{H^d, H^d} ^{-1}$. Then
$C_{U,V}=
(C_{V,U}) ^{-1}$,  for $ U, V = H $ or $H^d.$

\end {Lemma}
If $C_{H,H} = C_{H,H}^{-1}$, then we say that the braiding is symmetric on $H$.
If $C_{U,V} = C_{V,U}^{-1}$ for $U, V = H $ or $ H^d$, then we say that the braiding is symmetric on $H$ and $H^d$.
Throughout this section  we assume that the braiding is
symmetric on $H$ and $H^d$.

\begin {Lemma} \label {1.3}
(i)  $(H^{d}, \rightharpoonup )$  is a left $H$-module
algebra under the module operation.
 $\rightharpoonup \ = (H^d
\otimes <, >)(H^d \otimes C) (C
\otimes H^d  ) (H \otimes \Delta )$.

(ii)  $(H, \rightharpoonup )$  is a left  $H^{ d}$-module
algebra under the module operation
$\rightharpoonup \  = (H
\otimes <, >) (C \otimes H ) (H^d \otimes \Delta )$.
\end {Lemma}
Consequently, we can construct two smssh products $H\#H^d$ and $H^d\#H.$

\begin {Definition} \label  {1.4} Assume the category ${\cal C}$ is a subcategory of category ${\cal D }$.
We say that $CRL$-condition holds on $H$ and $H^d$ under $<, >$ if  the following conditions are satisfied:

(i) $E =: End _{\cal D} \ H$ is an algebra under multiplication of  composition  in ${\cal C}$ and  there exists a morphism $val: E \otimes H \rightarrow H$ in
${\cal D}$ such that
$val (f \otimes H) = val ( g \otimes H ) $ implies $f= g$ for any two morphisms $f$ and $g$ in ${\cal C}$ from  $U$ to $E$,
where $U$ is an object in ${\cal C}.$

(ii) There are two morphisms
 $\rho : H ^d \# H \rightarrow E$  and $\lambda :  H  \# H ^d \rightarrow E$
in ${\cal D}$ such that
$val (\lambda  \otimes H)=  (m\otimes <, > )(H \otimes C \otimes H )(H \otimes H^d \otimes \Delta )$ and
$val (\rho \otimes H)=  m (<, > \otimes C )(H^d \otimes C \otimes H)     (H \otimes H^d \otimes \Delta ).$

(iii) $Im ( \lambda )$ is an object in ${\cal D}$ and there exists a morphism $\bar \lambda $  in ${\cal D}$ from $Im ( \lambda )$ to $H\# H^d$ such
that $\bar \lambda \lambda = id _{H \# H^d}. $

(iv)
 $C_{E, V} (\lambda  \otimes V)=  (V \otimes \lambda )C _{H\# H^d , V}$ and
$ C_{E, V} (\rho  \otimes V)=  (V  \otimes \rho )C _{H^d\# H , V}$ for any object $V$ in ${\cal C}.$

\end {Definition}



\begin {Lemma}  \label {1.5}
$\lambda $ is an algebra morphism from $H\#H
^{d}$ to $E$ and $\rho $ is an
anti-algebra morphism from $H ^d\#H$ to $E$.
\end {Lemma}
{\bf Proof.}
We only need show that
$ val ((\lambda  \otimes H) (m \otimes H))= val (( m \otimes H) (\lambda \otimes \lambda \otimes H))$ and
$ val ((\rho  \otimes H) (m \otimes H))= val (( m \otimes H) (\rho  \otimes \rho \otimes H)(C \otimes H)).
$ The proof is similar to that of  \cite [Lemma 9.4.2]{Mo93}. $\Box$

\begin {Lemma} \label {1.6} The following relation  holds:
 $ m (\lambda \otimes \rho ) =  m (\rho \otimes \lambda )(H^d \otimes
\rightharpoonup \otimes \leftharpoonup \otimes H^d )(H^d \otimes
C_{H, H^d} \otimes C_{H^d,H} \otimes H^d) (H^d \otimes
H \otimes C_{H^d, H^d} \otimes H \otimes H^d )(H^d
\otimes C_{H^d,H} \otimes C_{H, H^d} \otimes H^d)(
 C_{H^d,H^d}\otimes  H \otimes H \otimes C_{H^d, H^d}
 )(S \otimes H^d  \otimes H  \otimes H  \otimes H^d  \otimes H^d )(\Delta  \otimes H
  \otimes H  \otimes \Delta  ) (H^d \otimes C_{H,H}\otimes H ) (C_{H, H^d } \otimes
  C_{ H^d,H })(H \otimes C_{H ^d,H ^d}\otimes H  ) $.
\end {Lemma}

{\bf Proof.} We show the relation after the following five steps. First we  check that the relation holds
on $(H\otimes \eta _{H^d}) \otimes ( \eta _{H^d}\otimes H)$,
$(\eta _H\otimes H^d)\otimes ( H^d\otimes \eta _H)$,
$(H\otimes \eta _{H^d})\otimes ( H^d\otimes \eta _H)$ and  $(\eta _H\otimes H^d )\otimes ( \eta _{H^d} \otimes H)$, respectively.
Using these
we check that the relation holds on $(H\otimes H^d)\otimes ( H^d\otimes H)$. $\Box$

\begin {Lemma} \label {1.7}
$R \# H$ becomes    an  $H^{ d}$-module algebra under the
 module operation  $\rightharpoonup ' \ =  (R \otimes
\rightharpoonup )(C_{H^{d}, R} \otimes H ).$
\end {Lemma}

{\bf Proof.} It is straightforward. $\Box$

Consequently, we obtain another  smash product $(R \# H)\# H^d.$

If $(R , \psi )$ is a right $H^{d}$-comodule algebra, then $(R, \alpha )$
becomes a left $H$-module algebra
(see \cite [Lemma 1.6.4] {Ma95b}) under the module operation:
$\alpha =(R \otimes  <, > )(R \otimes C_{H, H^d})(C_{H, R} \otimes R )(H \otimes \psi ) .$
          \begin {Theorem} \label {1.8}
 Let  $H$ be a   Hopf algebra. Assume that the $CRL$-condition holds on $H$ and $H^d$ under $<,>$,
and both $H$ and $H^{d}$ have invertible antipodes. Let $R$
be an $H^{d}$-comodule algebra, so that  $R$ is an $H$-module algebra defined as above. Let $H^{d}$ act on $R\#H$
by acting trivially on $R$ and via $\rightharpoonup $ on $H$,
then
  $$   (R \# H)\# H^{d}   \cong R \otimes (H
  \# H^{d} ) \hbox { \ \ \  as }  \hbox { algebras in }  {\cal D}.$$
In addition, if $\bar \lambda \rho (id _{H^d} \otimes \eta _H)$ is a morphism in ${\cal C}$ from $H^d$ to $H\# H^d,$
then the above isomorphism is one in ${\cal C}.$

\end {Theorem}
{\bf Proof. }
By (CRL)-condition,  there exists a morphism $\bar \lambda $ in ${\cal D}$
from $Im (\lambda )$ to
$H\# H^{d}$ such that  $\bar \lambda \lambda = id _{H \# H^{d}}$.
 We first define a morphism  $w=\bar  \lambda \rho (S^{-1}
\otimes \eta _H )$ from $H^{d} $
to $H \# H^{d}$. Since $\rho $ and $S^{-1}$ are
anti-algebra morphisms by Lemma \ref {1.5}, $w$ is an algebra morphism.

We now define two morphisms  $ \Phi = (R \otimes m _{H \# H^d})(R\otimes  w\otimes H \otimes H^d )
(\psi \otimes H \otimes H^d )$ from
$(R\#H)\# H^d$ to $R \otimes (H \# H^d)$ and
$ \Psi = (R \otimes m _{H \# H^d})(R\otimes  w\otimes H \otimes H^d )(R \otimes S \otimes H \otimes H^d ) (\psi \otimes H \otimes H^d )$
 from $R \otimes (H \#
H^{d})$ to $(R\#H)\# H^d.$
 It is straightforward to verify  that $\Phi \Psi = id  $ and $\Psi \Phi =
id$.
To see that $\Phi$ is an algebra morphism, we only need to show
that $\Phi ' = (R \otimes \lambda )\Phi $ is an algebra morphism.
Set  $\xi = (R \otimes \rho  )(R \otimes S^{-1} \otimes \eta _H) \psi $ from $R$ to $R \otimes (H \# \otimes H^d )$.
We have that $\xi $ is an algebra morphism
and
 $\Phi ' = (R \otimes m) (\xi \otimes \lambda ).$ Using Lemma
\ref {1.6}, we can show
  \begin {eqnarray*}
 &{\ }&(R \otimes m)(C\otimes E ) (\lambda \otimes \xi ) = \\
 &{\ }&(R  \otimes  m)
( \xi \otimes\lambda ) (\alpha \otimes H \otimes H^d  ) (H
\otimes C_{H, R} \otimes H^d )(\Delta  \otimes C_{H^d,
R} ) \ \ \ \ \ \ (*).
\end {eqnarray*}
 We now show  that
$\Phi'$ is an algebra morphism. Indeed,
 \begin {eqnarray*}
&{\ }& m _{R \otimes E}(\Phi ' \otimes \Phi ' )\\
 &=&( m \otimes m)(R \otimes C_{E, R} \otimes E )(m \otimes m)(\xi \otimes \lambda \otimes \xi
\otimes \lambda )\\
&=& (R \otimes m) ( m \otimes E \otimes  m)(R \otimes C_{E, R} \otimes m \otimes E)
(R \otimes E \otimes C_{E, R} \otimes E \otimes E)\\
&{\ }&(\xi \otimes \lambda \otimes \xi
\otimes \lambda )
(R \otimes m)(m \otimes E \otimes  m)(R \otimes C_{E, R} \otimes m \otimes E) \\
&{\ }&( R \otimes E \otimes  \xi \otimes\lambda \otimes E)
 (R \otimes E \otimes \alpha \otimes H \otimes H^d \otimes E  )\\
&=& ( R \otimes E \otimes H
\otimes C_{H, R} \otimes H^d \otimes E)
(\xi \otimes \Delta \otimes C_{ H^d ,  R } \otimes \lambda ) \hbox { \ \ by (*)}\\
&=& (m \otimes m)(R \otimes C_{E,R} \otimes m)( R \otimes E\otimes  \xi \otimes \lambda )
(R \otimes E \otimes \alpha \otimes m_{H \#H^d}  )\\
&{\ }& ( R \otimes E \otimes H
\otimes C_{H, R} \otimes H^d \otimes H \otimes H^d)\\
&{\ }& (\xi \otimes \Delta \otimes C_{H^d, R} \otimes H \otimes H^d ) \hbox { \ \ (by Lemma \ref {1.5})}\\
&=& (R \otimes m)(\xi  \otimes \lambda )( m   \otimes H\otimes H^d)
(R  \otimes \alpha \otimes m_{H \#H^d}  )
 ( R  \otimes H
\otimes C_{H, R} \otimes H^d \otimes H \otimes H^d)\\
&{\ }& (R \otimes \Delta \otimes C _{H^d , R} \otimes H \otimes H^d ) \hbox { \ \ (since
} \xi  \hbox { is algebraic ) }\\
&=&\Phi ' m_{(R \# H) \# H^d}.
 \end {eqnarray*}
Thus $\Phi '$ is algebraic and $\Phi $ is also algebraic.
If $\bar \lambda \rho (id _{H^d} \otimes \eta _H)$ is a morphism in ${\cal C},$
then $\Phi $ is an isomorphism  in ${\cal C}.$
$\Box$

We obtain the following by  Theorem \ref {1.8}.

  \begin {Corollary} \label {1.9}
 Let  $H$ be a finite braided  Hopf algebra with a left dual $H^*$.
If the braiding is symmetric on $H$, then
  $$   (R \# H)\# H^{*}   \cong R \otimes (H   \# H^{*} ) \hbox { \ \ \  as }  \hbox { algebras  in } {\cal C}. $$
\end {Corollary}
This corollary reproduces   the main result in \cite {Zh03}.

\section {Duality theorems in the Yetter-Drinfeld module category} \label {e2}

In this section, we present the duality theorem for braided Hopf algebras
in the Yetter-Drinfeld module category $(^B_B {\cal YD}, C)$.
Throughout this section, $H$ is a braided Hopf algebra in $(^B_B{\cal YD}, C)$ with finite-dimensional Hopf algebra $B$
and $H^d$ is a quasi-dual of $H$ under a left faithful $<, >$ (i.e. $<x, H > = 0$ implies $x =0$) such that
$<b\cdot f, x> = <f, S(b)\cdot x>$ and $\sum <f _{(0)}, x>f_{(-1)} = \sum <f, x_{(0)}> S^{-1}(x_{(-1)}) $ for any $x\in H,
b \in B, f \in H^d.$
 Let $b_B$ denote the coevaluation of $B$ and   $<, >_{ev}$  the ordinary evaluation of any spaces.

\begin {Lemma} \label {2.1}
(i)  If $(V, \alpha_V, \phi _V)$ and $(W, \alpha_W, \phi _W)$ are two Yetter-Drinfeld modules over $B$, then
 $Hom _k (V, W)$ is a Yetter-Drinfeld module under the following module operation and comodule operation:
 $(b \cdot f)(x) = \sum b_2\cdot f(S(b_1)\cdot x )$ and
$\phi (f) = (S^{-1} \otimes \hat \alpha ) (b_B \otimes f )$,
 where  $ \hat \alpha $ is defined by
$(b^* \cdot f)(x) = <b^*, x_{(-1)} S(f (x_{(0)})_{(-1)})>_{ev}(f (x_{(0)}))_{(0)}$ for any $x\in V,
f \in Hom _k (V,W), b^* \in B^*.$ In particular, if $V$ is an object in $^B_B {\cal YD}$, then so is $V^*$.

(ii) If $V$ is an object  in $(^B_B{\cal YD}, C)$, then  $V^*$ is object   in $(^B_B{\cal YD}, C)$ and the evaluation
$<,>_{ev}$ is a morphism in $(^B_B{\cal YD}, C).$

(iii) If the braiding is symmetric on $V$, then it is symmetric on $V$ and $V^*$.

 \end {Lemma}

{\bf Proof.}
 (i)  It is clear  that
$\sum f _{(-1)} f _{(0)}(x) = \sum (f (x_{(0)}))_{(-1)} S^{-1}(x_{(-1)}) \otimes (f (x_{(0)}))_{(0)} $
for any $x\in V, f \in Hom _k (V, W), b \in B.$  Using this, we can show that $Hom _k (V, W)$ is a $B$-comodule.
Similarly, we can show that $Hom _k (V, W)$ is a $B$-module. We now show that
$$ \sum (b \cdot f )_{(-1)} \otimes (b \cdot f )_{(0)}= b_1 f _{(-1)}S(b_3) \otimes b_2 \cdot f _{(0)} \ \ \ \ \ \ \ \ (*)$$
for any $f \in Hom _k(V, W), b \in B$. For any $x\in V,$  see that
 \begin {eqnarray*}
 \sum (b \cdot f )_{(-1)} \otimes (b \cdot f )_{(0)}(x)&=& \sum b_1 (f (S(b_4))\cdot x_{(0)})_{(-1)} S(b_3)S^{-1}(x_{(-1)})\\
&{\ }& \otimes b_2 \cdot (f (S(b_4))\cdot x_{(0})_{(0)}   \hbox {\ \ \ \ \ and } \\
b_1 f _{(-1)}S(b_3) \otimes (b_2 \cdot f _{(0)})(x)&=& b_1 f_{(-1)} S(b_4) \otimes  b_2 \cdot f_{(0)}
((S(b_3)\cdot x))\\
&=& \sum b _1 (f (S(b_4) \cdot x _{(0)}))_{(-1)} S(b_3)S^{-1}(x_{(-1)})b_5S(b_6) \\
&{\ }& \otimes
 b_2 \cdot (f (S(b_4) \cdot x _{(0)}))_{(0)} \\
&=& \sum b_1 (f (S(b_4))\cdot x_{(0)})_{(-1)} S(b_3)S^{-1}(x_{(-1)})\\
&{\ }& \otimes b_2 \cdot (f (S(b_4))\cdot x_{(0)})_{(0)}.
 \end {eqnarray*}
Thus (*) holds and $Hom _k (V, W)$ is a Yetter-Drinfeld module.

(ii) By (i), $V^*$ is a  Yetter- Drinfeld $B$-module.
Obviously, $<,>$ is a $B$-module homomorphism. In order to show that  $<,>$ is a $B$-comodule homomorphism, it is
enough to prove  that $\sum h_{(-1) }^* h _{(-1)} <h^*_{(0)}, h_{(0)}> = 1_B <h^*, h>$ for any $h^* \in V^*, h \in V.$
Indeed,  the left side $= \sum S^{-1}(h_{(-1) 2}) h _{(-1)1} <h^*, h_{(0)}> = 1 _B <h^*, h>.$ This
complete the proof.

(iii) It follows from Lemma \ref {1.2}.
 $\Box$

\begin {Lemma} \label {2.2}
Let $A$ be a braided algebra in $ {\cal C} =( {}^B_B{\cal YD}, C)$ and
$A^\circ _{\cal C}= \{f \in H^* \mid Ker (f) \hbox
{ contains an ideal
of finite codimension in } {}^B_B{\cal YD} \} $. Then  $ A^ { \circ } _{\cal C} $
is a braided coalgebra in $ (^B_B{\cal YD}, C)$, called the finite dual of $A$ in ${\cal C}$ and written as
$A^\circ$ in short. Moreover, if $H$
is a braided Hopf algebra in ${\cal C}$, then $H^\circ_{\cal C}$ is a braided Hopf algebra in ${\cal C}$.
\end {Lemma}

{\bf Proof.}  By Lemma \ref {2.1}, $A^*$ is a $B$-module and $B$-comodule. First we show that $A^\circ$ is an object in $^B_B{\cal YD}$. For any $f \in A^\circ$, there exists
an ideal $I$ of $A$ and $I$ is a $B$-submodule and a $B$-subcomodule of $A$ with finite codimension and $f(I)=0$. Since
$(b\cdot f) (x) = f (S(b)\cdot x) = 0$ for any $b \in B, x \in I$, we  have  $b\cdot f \in A^\circ $. Thus
$A^\circ$ is a $B$-submodule of $A^*$. By Lemma \ref {2.1}, we can assume $\phi _{A^*} (f) = \sum_i u_i \otimes v_i$ with linear independent $u_i's.$
Since $\sum_i u_i v_i (x) = \sum f (x_{(0)})S^{-1} (x_{(-1)})=0 $  for any $x\in I$, we have that $v_i (x)=0$ and
$v_i(I)=0$, which implies $v_i \in A^\circ$. thus $A^{\circ}$ is a $B$-subcomodule of $A^*.$

We next  show that $A^\circ \otimes A^\circ = (A\otimes A)^\circ $ and $m^* (A^\circ) \subseteq A^\circ \otimes
A^\circ$ by using the method  similar to the proof in \cite   [Lemma 1.5.2] {DNR01}.
To show that $m^*, \eta ^*$ are morphisms in $^B_B {\cal YD},$ we only need show that
if $f $ is a morphism from $U$ to $V$ in $ {}^B_B{\cal YD}$, then $f^*$ is a morphism from $V^*$ to $U^*$ in
$ {}^B_B{\cal YD}$. Indeed, for any $v^* \in V^*, u \in U, b\in B,$ see that
\begin {eqnarray*}
((b \cdot f ^* (v^*))(u) &=& (f^* (v^*))(S(b) \cdot u)\\
&=& v^* (f (S(b) \cdot u)) \\
&=& (f^* (b \cdot v^*))(u).  \end {eqnarray*}
Thus $(b \cdot f^*) (v^*)= f^* (b\cdot v^*)$ and $f^*$ is a $B$- module homomorphism. Similarly, we can show that
$f^*$ is a $B$-comodule homomorphism. Consequently, $(A^\circ , m^*, \eta ^*)$ is a braided coalgebra in $(^B_B{\cal YD}, C).$
Finally we can similarly complete the other.
$\Box$

Let $\lambda '$ denote the $k$-linaer map from $H\# H^{d}$ to $End _k H$ by sending $h \otimes h^{d}$ to $\lambda ' (h \# h^{d})$
with $\lambda ' (h \# h^{d})(x)= h <h^{d}, x>$ for any $x\in H, h\in H, h^{d} \in H^{d}.$
Obviously, $\lambda '$ is  an injective $k$-linear map, so we can view  $H\# H^d$ as a subspace of $End _k H.$
Now we define $\lambda $
and $\rho $. For any $h, x \in H, f \in H^d$, $(\lambda (h \# f))(x) =:
(m \otimes <,> )(H \otimes C \otimes H ) (H \otimes H^d \otimes \Delta ) (h \otimes f \otimes x)
= \sum <f,x_{2(0)}>h (S^{-1}(x_{2(-1)} \cdot x_1))$ and $(\rho (f \# h))(x) =:
(<,> \otimes m )(H^d \otimes H \otimes C  )(H ^d \otimes C \otimes H ) (H^d \otimes H \otimes \Delta ) (h \otimes f \otimes x)
= \sum <f, h _{(-1) 1} \cdot x_1>( h_{(-1)2} \cdot x_2 ) h_{(0)}.$

Let  ${\cal D}$ denote  the category of  vector spaces and  ${\cal
C} = {}^B_B {\cal YD}$. Define $val (f \otimes x) =f(x)$ for any
$f \in E, x \in H.$ If $\rho (H^d\# 1) \subseteq \lambda (H\#
H^d)$ then  we say that $RL$-condition holds on $H$ and $H^d$
under $<,>$.

\begin {Lemma} \label {2.3} Let $H$ be a braided Hopf algebra in $(^B_B{\cal YD}, C)$ with $C_{H,H} = C_{H,H}^{-1}.$

(i)  If the antipode of $H$ is  invertible, then there exists $k$-linear map  $\bar \lambda $
 from $Im \lambda $ to $H\# H^{d}$ such that $\bar \lambda \lambda  = id _{H\# H^{d}}$.

(ii)  If $H$ is quantum  cocommutative, then $RL$-condition holds on $H$ and $H^d$ under $<,>$.

(iii) $C_{E, V} (\lambda  \otimes V)=  (V \otimes \lambda )C _{H\# H^d , V}$ and
$ C_{E, V} (\rho  \otimes V)=  (V  \otimes \rho )C _{H^d\# H , V}$ for any object $V$ in ${\cal C}.$

(iv) $E = End _k H$ is an algebra in $^B_B {\cal YD}.$

(v) If $B$ is a commutative and  cocommutative finite-dimensional Hopf algebra and $H$ has an invertible antipode, then
$\bar \lambda \rho (id _{H^d} \otimes \eta _H)$ is a morphism in $^B_B {\cal YD}.$
\end {Lemma}

{\bf Proof.} (i)
We  define a $k$-linear map $ \bar \lambda $ from $Im \lambda $ to $H\# H^d$ as follows:
$\bar \lambda  (f) (x) = m (f \otimes H) C (S^{-1} \otimes H)\Delta (x)$  for any $f \in Im \lambda , x \in H.$
We can show that $\bar \lambda  \lambda  = id _{H\# H^d}.$  Indeed,
for any $h\in H, h^{d} \in H^{d}, x \in H$, we have
\begin {eqnarray*}
\bar  \lambda  \lambda (h \# h^{d}) (x) &=& m(m \otimes <,> \otimes H)
(H \otimes C  \otimes H \otimes H)(H  \otimes H^{d} \otimes \Delta \otimes H)\\
&{\ }&(H  \otimes H^{d} \otimes C)
(H  \otimes H^{d} \otimes S^{-1}\otimes H)(H  \otimes H^{d} \otimes \Delta ) (h \# h^{d}) (x)\\
&=& \lambda ' (h \# h^{d}) (x).
\end {eqnarray*}
Thus $ \bar \lambda \lambda  = id _{H\# H^{d}}.$

(ii) It follows from the simple fact
$\rho (f \# 1) = \lambda (1 \# f)$ for any $f \in H^d$.

(iii) We only show that $C_{E, V} (\lambda  \otimes V)=  (V \otimes \lambda )C _{H\# H^d , V}.$ Indeed,
 for any
$h, x \in H, v\in V$, see that
\begin {eqnarray*}
(V \otimes &val& ) (C_{E, V}\otimes H) (\lambda  \otimes V \otimes H) (h \otimes f \otimes \otimes v \otimes  x)\\
&=& <f , x_{(0)2(0)}>( h_{(-1)}S^{-1}(x _{(0)2(-1)})_1 x_{(0)1(-1)} S^{-1}(x _ {(0)1(-1)}S( S^{-1}(x_{(0)2(-1)}))_3) \cdot v \\
&{\ }& \otimes h _{(0)} S^{-1}(x_{(0)2(-1)})_2 \cdot x_{(0)1(0)}) \hbox {\ \ }\\
&=&<f , x_{2(0)}> h_{(-1)} S^{-1}(x _{2 (-1)4})x_{1(-1)2}\underline {x_{2(-1)2} S^{-1}(x_{2(-1)1})}
 S^{-1} (x_{1 (-1)1}) \cdot v \\
&{\ }& \otimes h_{(0)} (x _{2 (-1)3}\cdot x_{1(0)}) \\
&=& <f , x_{2(0)}> h_{(-1)} S^{-1}(x _{1 (-1)2}) \cdot v \otimes h_{(0)}(S^{-1}(x_{2(-1)1}) \cdot x_1) \ \ \hbox { and }\\
 (V \otimes &val& ) (V \otimes \lambda \otimes H)(C _{H\# H^d , V} \otimes H)(h \otimes f \otimes \otimes v \otimes  x)\\
&=& <f_{(0)}, x_{2(0)}>h_{(-1)}\cdot (f_{(-1)} \cdot v) \otimes h _{(0)} (S^{-1}(x_{2(-1)})\cdot x_1 ) \\
&=& <f , x_{2(0)}> h_{(-1)} S^{-1}(x _{1 (-1)2}) \cdot v \otimes h_{(0)}(S^{-1}(x_{2(-1)1}) \cdot x_1).
\end {eqnarray*}
Thus $C_{E, V} (\lambda  \otimes V)=  (V \otimes \lambda )C _{H\# H^d , V}.$
$\Box$

(iv) It is straightforward.

(v) Let $\mu $ denote  $ \bar \lambda \rho (id _{H^d} \otimes \eta _H).$ We only show that $\mu$ is a $B$-module
homomorphism.
For any $x \in H, f \in H^d, b \in B,$ since $B$ is commutative and cocommutative, we have
$\sum (b\cdot x)_{(-1)} \otimes (b\cdot x)_{(0)} = \sum x_{(-1)} \otimes (b\cdot x_{(0)})$  and

\begin {eqnarray*}
(\mu (b \cdot f ))(x) &=& \sum <f ,S(b) x_{1(-1)1} \cdot x_2 > (x_{1(-1)2} \cdot x_3)S(x_{1 (0)})\\
(b \cdot \mu (f))(x) &=&    <f , x _{1 (-1)1}S(b _4) \cdot x_2  > (b_1 x_{1 (-1)2} S(b_5) \cdot x_3)(b_2S(b_3)\cdot
 S(x_{1(0)}))\\
&=& \sum <f ,S(b) x_{1(-1)1} \cdot x_2 > (x_{1(-1)2} \cdot x_3)S(x_{1 (0)}).\\
\end {eqnarray*}
Thus $\mu$ is $B$-module homomorphism. $\Box$

Every $B$-module  category  $(_B{\cal M}, C^R)$ determined by quasitriangulr Hopf algebra $(B,R)$
 is a full subcategory  of  the Yetter-Drinfeld
module category  $(^B_B{\cal YD}, C)$.   Indeed, for any $B$-module $(V, \alpha   )$, define $\phi (v) = \sum
R_i ^{(2)} \otimes R^{(1)}_i \cdot v$  for any $v \in V$, where $R = \sum _i R_i ^{(1)} \otimes R_i ^{(2)}$.
It is easy to check that $(V, \alpha , \phi )$ is a Yetter-Drinfeld $B$-module. Similarly,
every $B$-comodule category $(^B{\cal M}, C^r)$ determined by coquasitriangulr Hopf algebra $(B,r)$
 is a full subcategory of the  Yetter-Drinfeld
module category $(^B_B{\cal YD}, C)$.


  \begin {Theorem} \label {2.4}
 Let  $H$ be a  braided  Hopf algebra in
$(^B_B{\cal YD}, C)$ with finite-dimensional  $B$ and $C_{H,H}= C_{H,H} ^{-1}$. Assume that  $RL$-condition holds on $H$ and
$H^d$ under left faithful $ <, >$, and both $H$ and $H^{d}$ have invertible antipodes. Let $R$
be an $H^{d}$-comodule algebra, so that  $R$ is an $H$-module algebra defined as above. Let $H^{d}$ act on $R\#H$
by acting trivially on $R$ and via $\rightharpoonup $ on $H$. Then
  $$   (R \# H)\# H^{d}   \cong R \otimes (H   \# H^{d} ) \hbox { \ \ \  as } k\hbox {-algebras.   } $$
Moreover, if $B$ is commutative and  cocommutative, then the above isomorphism is one as algebras in $^B_B {\cal YD}.$
\end {Theorem}

{\bf Proof.}
It follows from Lemma \ref {2.1} and Lemma \ref {2.3} that (CRL)-condition in Definition \ref {1.4}
is satisfied. Considering Theorem \ref {1.8}, we complete the proof. $\Box$

  \begin {Corollary} (Duality Theorem) \label {2.5}
 Let  $H$ be a  braided  Hopf algebra in
$(^B_B{\cal YD}, C)$ with finite-dimensional  $B$ and $C_{H, H} = C_{H,H}^{-1}$.
Assume that $U$ is a braided Hopf subalgebra of $H^\circ  $
and $RL$-condition holds on $H$ and $U$ under evaluation $<, >_{ev}$, and  $H$ has  invertible antipode. Let $R$
be an $U$-comodule algebra, so that  $R$ is an $H$-module algebra defined as above. Let $U$ act on $R\#H$
by acting trivially on $R$ and via $\rightharpoonup $ on $H$,
then
  $$   (R \# H)\# U   \cong R \otimes (H   \# U ) \hbox { \ \ \  as }  k\hbox {-algebras. }   $$
Moreover, if $B$ is commutative and  cocommutative, then the above isomorphism is one as algebras in $^B_B {\cal YD}.$
\end {Corollary}

{\bf Proof.} It is clear that $U$ is a quasi-dual of $H$ under evaluation $<,> _{ev} $, which is a left faithful.
$U$ has an invertible antipode since $H$ has an invertible antipode.
By Theorem \ref {2.4},  we  complete the proof. $\Box$

  \begin {Corollary} ( Second Duality Theorem )\label {2.6}
 Let  $H$ be a  braided  Hopf algebra in
$(^B_B{\cal YD}, C)$ with finite-dimensional  $B$ and $C_{H, H} = C_{H,H}^{-1}$. Let $U$ and $V$ be
a braided Hopf subalgebra
of $H^\circ$  and  $H$ with invertible antipodes, respectively.
Assume that $RL$-condition holds on $U$ and $V$ under $<, > = <, >_{ev} C$, and $U$ is  dense in $H^*$.
 Let $R$
be an $V$-comodule algebra, so that  $R$ is an $U$-module algebra defined as above. Let $V$ act on $R\# U$
by acting trivially on $R$ and via $\rightharpoonup $ on $U$,
then
  $$   (R \# U)\# V   \cong R \otimes (U   \# V ) \hbox { \ \ \  as } k\hbox {-algebras.  } $$
Moreover, if $B$ is commutative and  cocommutative, then the above isomorphism is one as algebras in $^B_B {\cal YD}.$
\end {Corollary}

{\bf Proof.} It is clear that $V$ is a quasi-dual of $U$ under $<, > = <,> _{ev}C_{V, U} $ and $<,> $
is left faithful since $U$ is dense. By Theorem \ref {2.4},we can complete the proof. $\Box$


{\bf Remark:} In Corollary \ref {2.5} $H$  can be replaced by  Hopf subalgebra $V$ of $H$ when $U \subseteq V^*$ and
$V$ has an
invertible antipode.


\begin {Example} \label {2.7} Let $H$ be  a quantum cocommutative braided Hopf algebra in ${\cal C}={}^B_B{\cal YD}$
with  finite-dimensional commutative and cocomutative $B$ and $C_{H,H}= C_{H,H}^{-1}$ (for example, $H$ is the
unversal enveloping algebra of a Lie superalgebra). Set $U =H^\circ _{\cal C} = A$. It is clear that $(A, \phi )$ is a right  $U$-comodule algebra with $\phi = \Delta $.
By Lemma \ref {2.3}, The $RL$-condition holds on $H$ and $U$ under evaluation $<,>_{ev}$. By
Corollary \ref {2.5}, we have
 $$   (R \# U)\# H   \cong R \otimes (U   \# H ) \hbox { \ \ \  as } \hbox {algebras in } {}^B_B{\cal YD}.  $$
\end {Example}




{\bf Remark:} Although we have an efficient Sweedler's method (see
 \cite {Sw69}) to present co-operations,
 we suggest that readers use braiding diagrams to check all of our proofs because they are clearer.

\vskip 0.5cm
{\bf Acknowledgement }: The work is supported by Australian Research Council.
S.C.Z. thanks the Department of Mathematics, University of
Queensland for hospitality.

\end {document}